# Aletheia, double negation and negation

Antonino Drago
University "Federico II" Naples – I-, drago@unina.it

**Abstract.** The main result of present paper is that the definition of negation has to be referred to the totality of a theory and at last to what is defined as the organization of a scientific theory; in other words, negation's definition is of a structural kind, rather than of an objective kind or a subjective kind. The paper starts by remarking that the ancient Greek word for truth was "aletheia", which is a double negation, i.e. "unveiling". But, after Plato and Romans the affirmative meaning of the idea of truth prevailed. Not before the 1968 the double negation law was re-evaluated, since it was recognized that its failure represents more appropriately than the failure of of excluded middle law the borderline between classical logic and almost all non-classical kinds of logic. Moreover, the failure of this law is easily recognized within a (scientific) text; this fact allows a new kind of logical analysis of a text. As an example, the analysis of Kolmogorov's 1932 paper is summarized. It shows that he reasoned according to arguments of non-classical logic about the foundations of the intuitionist logic. Through a comparative analysis of all scientific theories which are based, like the previous one, on a general problem my previous paper suggested an alternative model of organization of a theory to the deductive-axiomatic one. Then general informal logic is founded according to this model. The negation is defined as a unary operation which, under the problem of deciding whether a doubly negated proposition is equal to the corresponding affirmative proposition or not, leads to a subdivision into classical logic and intuitionist logic. Hence, negation takes two specific and mutually incompatible meanings; in classical logic it leads to the non-contradiction principle, whereas in intuitionist logic it does not opposes to the corresponding affirmation. This essentially ambiguous nature of the negation gives reason of the insufficiency of Natural Deduction.

### 1. A dichotomy on the kind of logic: *Aletheia* vs. truth

Let us consider how in the past times the notion of "truth" has been expressed by a word. The ancient Greeks had the word *alethéia* and the Romans the word *veritas*. The first word comes from a privative α and λανθάνω which means "I am hidden" (like the mythological river *Lete*, ἡ Λήθη, of which Plato speaks for example in the tenth book of *Republic*, is the river of oblivion). So ἀλήθεια is the state of "not being hidden"; we could say "unveiling". Thus, it has to be regarded not as a simple factual reality or as a hypostatized idea, but as a dynamic act, through which the recognition and the refutation of errors occur: hence, not a static and defined thought once and for all, but a revelatory movement. One easily recognizes that natural language includes several doubly negated words (e.g. in-nocent, in-nocuous, non-violent, etc.).

The term *veritas*, on the other hand, comes from the Balkan and Slavic areas and, according to its original meaning, means "faith" (the fact that the wedding ring in Italian is called "vera" confirms this ethimology) and therefore the fact is true without any need for critical scrutiny.

Plato changed the meaning of aletheia into a fixed, metaphysical meaning. Aristotle was immersed within the new meaning; he stated: "To say of what is that it is, and of what is not that it is not, is true." (Metaphysics 1011b) This fact means that the history of the Greek-Roman civilization was marked by an unrecognized inner conflict between two meanings of the same word, truth. This conflict was won by the idealistic meaning of truth.

Owing to the introduction of an ever great number of abstract ideas, the natural languages so much lost a strict correspondence with facts, that often no attention is done whether a DN proposition is supported by evidence or not; this trend was not shared by Court's language which has to be relied upon safe evidence and moreover the language of experimental science which is essentially linked to facts (although also natural science and mathematics introduced some idealistic ideas).

Indeed, Galilei's revolution in the methodology of science required that the truth of a scientific proposition be supported by experimental data. Yet, this requirement has been ignored by the language of philosophers and in partial way also scientists, because since the birth of modern Mathematics, mathematicians made use of metaphysical infinitesimals (and then idealistic axioms, as Zermelo's) and subsequently also the logicians made a unrestricted use of predicate calculus. Hence, linguists had good reasons for resisting to this requirement.

However, in past century some Intuitionists (Dummett 1977) and Constructivists (Bishop 1967, pp. 1-10) scholars claimed inside respectively Logic and Mathematics a similar requirement to the Galilean one; when a proposition does not satisfy such a requirement, its doubly negated versions can be however stated. Moreover, Dag Prawitz and Peer Melmnaas (1968) stated that the validity or not of the law of double negation separates classical logic (CL) from almost all non-classical kinds of logic, first of all the intuitionist one, and then the minimal one, the positive one, etc.

In the light of this result the two old words of truth also represent e formally logical difference, because the double negation 'aletheia' belongs to non-classicl logic (NCL) and 'truth' to CL. The two words represent a radical variation in the logical meaning of a same notion shared by the two logical theories.

In retrospect, we see that the culture of Western civilization was linked to the word *veritas* and hence to only CL (recall Immanuel Kant's faith on CL as already completed since longtime[1]).

Truly, Thomas Aquinates had suggested a new definition: "*Veritas est adaequatio rei et intellectus*". It leaves to intend that truth is reachable by means of an approximation process, which at last achieves an affirmative term, as "equality". Instead, *adaequatio* is a modal word; it means "to become similar". Being the modal logic translatable into intuitionist logic (IL) through the S4 model (Hughes and Cresswell 1996, pp. 224ff.), this novelty actually generated a compromise between CL and NCL, which represented an ambiguous way of intending the truth. No surprise if the debate on the notion of truth, although absorbing the reflections of the major minds of all past centuries, did not achieved a stable and commonly accepted result and many troubles arise also at present time from the confusion of these two conceptions of truth.[2] No surprise if the notion of negation represents a great problem.

**2. Analysis of a scientific text through the law of the double negation: How many double negations?**

The shift from the law of the excluded middle (LEM) to the law of double negation (DNL) as the borderline between CL and the main kinds of NCL is of a great importance.

A very important fact is that within a text LEM is rarely manifested (possibly as dilemmas); hence, almost never a text reveals through the LEM whether the underlying logic is different from the classical one. Instead, the reading of a text easily manifests author's use of a NCL by means of the failure of the DNL, provided that one recognizes two negations within a proposition and moreover proves that its meaning is different from that of the corresponding affirmative proposition (DNP). For instance, Hyppocrates' maxim: "First, do not harm" (≠ "Do the good"); or the Decalog's maxim "Thou do not kill" (≠ "Save life").

From my experience of thirty years of recognizing double negations in the original texts of logic, mathematics, physics and chemistry, I extracted the following rules:

1) One has to exclude the rhetoric use of the double negations (e.g. "I have nothing else than 10$", which corresponds to the proposition obtained by making use of the law classical of double negation: "I have 10$").

2) The word "only", i.e. nothing else, works as a DNP.

---

[1] Truly, Hegel tried to change the kind of logic. Unfortunately his result was so confused that it is recalled as a negative instance. Instead Cusanus in 15th Century made consistent use of non-classical logic (Drago 2010)

[2] A longstanding tradition of English linguists suggests what Lawrence Horn called a "dogma" about double negations. (Horn 2001, pp. 79ff.; Horn 2008).This linguistic dogma asserts the absolute validity of the double negation law: whenever a DNP is found in a text, it has to be changed into the corresponding affirmative proposition, because those who speak by means of DNPs want to be, for instance, unclear. Evidence for this "dogma" is the small number of studies on double negations in comparison with the innumerable studies on the single negation.

3) Modal words "must", "possible", "necessary", etc. work as DNPs, owing to the S4 translation between modal logic and IL.

4) A comparison as "no more... than...", "no less... than..." works as a DNP, because it is not equivalent to the identity (in mathematical terms it corresponds to the Stone representation of IL by means of open sets).

5) An interrogative negative proposition is a DNP, when the answer is implicitly intended as a negative one.

6) A proposition may understood the second negation; it is recognized through the meaning of the entire discourse: Popper: "Science is fallible [owing to negative experimental data]"; Court's judgment: "Acquitted for insufficient evidence [of guilt]"; Jonas: "The ethics of the fear [of the mankind's suicide]"; hence, one has to grasp the full meaning of a negated proposition in order to discover a hidden second negation.

The character of being a negative word may depend on the context; for example in science the words "change", "different", "variation", etc. all work as negative words because each of them represents a problem which requires to be explained by a scientist.

Through the DNL one easily recognizes in an original text founding a scientific theory the kind of logic that its author is following; one has to inspect whether there exist relevant DNPs. I discovered that there exist several scientific theories that have been presented by their respective authors through texts containing several DNPs. In other words, the presentations of several theories instantiated a specific use of non-classical logic, since they relied in an essential way on DNPs. This kind of analysis was successful in e.g. Avogadro's atomic theory, S. Carnot's thermodynamics, Lobachevsky's hyperbolic geometry, Einstein's 1905 paper on quanta, Kolmogorov's minimal logic and Markov's theory of constructive numbers, Computer science (Church's thesis).

### 3. A new model of organization of a scientific theory and the two basic dichotomies

I performed a systematic research on the original texts of past theories of the hard science containing DNPs. These theories are the following ones: L. Carnot's theories of calculus and geometry, Lagrange's mechanics, Lavoisier's chemistry, Avogadro's atomic theory, S. Carnot's thermodynamics, Lobachevsky's non-Euclidean geometry, Galois' theory of groups, Klein's Erlanger program, Poincaré's theory of integer numbers, Einstein's theory of special relativity, both Planck's and Einstein's theories of quanta, Kolmogorov's foundation of minimal logie, Markov's theory of constructive functions. Remarkably, each of these theories is organized not in an axiomatic way, rather in order to solve a general problem (PO) on the background of the common knowledge shared about the field at issue. By ingenuity an author makes use of DNPs so that their mere sequence gives the logical thread of the entire illustration. This fact proves that The DNPs ply an essential role in the author's illustration.

Moreover, I discovered that an argumentation including a DNP, if correctly formulated, leads to an *ad absurdum* argument (AAA)[3]; as an instance S. Carnot's *ad absurdum* theorem (based on the DNP "A motion without an end is impossible" and proving a very concrete fact of machines!) is basic within thermodynamics, Here, NCL manifests an original way of arguing. An analysis of his original text shows that it includes seven AAAs. Also Lobachevsky's theory on the problem of parallel lines includes a similar chain of five AAAs.

The conclusion of such a chain of AAAs is a predicate - again a DNP – which constitutes a surmise for the resolution of the stated problem. Its translation in a logical formula gives in some theories $\neg \exists x \neg f(x)$, in the others $\neg \neg \forall x\, f(x)$. Remarkably, these two intuitionist predicates are equivalent, as Dummett's table (1977, p. 29) shows.

After this conclusion, it is easily recognised that the author of a PO theory refers to instead of the previous predicate, the corresponding affirmative, in order to deductively draw from it, according to CL, all the relevant consequences, to be then tested with reality. In other words, the author has assumed the corresponding affirmative predicate as a new hypothesis for a subsequent deductive system. The formula of this logical change is the following one:

$$\neg\neg \forall x\, f(x)\ (\leftrightarrow \neg \exists x \neg f(x))\ =>> \forall x\, f(x).$$

---

[3] Indirect proofs have been considered equivalent to classical ones through the use of DNL (see e.g. Gardiès 1991)

Through Dummett's table one easily shows that the above translation implies the translation of all intuitionist predicates (and propositions) on the same subject into the corresponding classical predicates (and propositions); hence, after this move an author correctly makes use of CL.

Markov's theory performs this translation (in order to obtain his well-known principle) by means of words which represent the principle of sufficient reason (PSR).(Markov 1962, p. 5) This principle was suggested by Leibniz in alternative to the non-contradiction principle, which governs an AO theory. By analysing Leibniz' text presenting the PSR, (Leibniz 1686) at glance a DNP is recognized "<u>Nothing</u> is <u>without</u> reason"; its logical formula is $\neg \exists X \neg f(X)$, i.e. the same as the predicate concluding the non-classical arguing in each PO theory, but now, owing to its nature of a principle, it is applied to whatsoever content.

Leibniz continues as follows: "..., or everything has its reason, although we are not always capable of discovering this reason. . ." Admirably, the last proposition cleverly explains why the affirmative proposition is not equivalent to the former proposition, a DNP. Owing to this inequivaence, this principle leads to enquire each case through a specific method. In formal terms, Leibniz' application of PSR to a predicate is represented by

$$\neg \exists X \neg f(X) =>> \forall X f(X).$$

This is the same translation of the final predicate of a PO theory. For the same above mentioned reason, this move translates all IL predicates, now on whatsoever argument, into CL predicates. Hence, it is the application of the PSR which governs the final step of a PO theory.

Markov suggests two constraints to the application of PSR to the conclusive predicate: i) to be the conclusion of an AAA and ii) to be decidable. Both are appropriate for the application of PSR to the final predicate of each PO scientific theory; indeed, the conclusive predicate results from (the last) AAA and moreover it is on operative - hence decidable – statement, owing to the scientific nature of the method solving the problem of a PO theory[4].

In sum, a formal model of a new organization of a scientific theory results[5]. It starts from a crucial problem whose resolution requires a new scientific method. It proceeds through DNPs, composing AAAs, linked together in a chain of derivations, concluding a doubly negated predicate, to which the author implicitly applies the methodological principle, the PSR, for obtaining a new hypothesis, from which to derive consequences to be tested with reality. Notice that theapplication of PSR constitutes the inverse (doubly) negative translation of CL into IL; a translation which is not obtained by the application of the Markov's principle. The above results are all illustrated by the papers (Drago 2012; Drago 2017a).

Notice that if the PSR is not recognized as a logical operation (as it was in past centuries), only CL is a formal logic, whereas IL, whose complete logical foundation requires the use of PSR, is no more than an intuitive, informal kind of logic; this was the attitude of logicians of the past centuries, at the cost to be blind to half part of logical arguing, depreciated as a merely intuitive, informal one (also Brouwer perceived it in this way). If, on the contrary PR is recognized as a logial operation, i.e. the translation of the IL into CL, then we obtain a pluralism of formalizations of general logic and hence a great empowerment of past logical arguing.

---

[4] Instead, Dummett's table shows that usual Markov's principle cannot translate the entire IL into CL.

[5] Unsuccessfully William Beth searched through a philosophical analysis, i.e. a non-deductive way of organizing a mathematical theory. But he ignored that as first D'Alembert suggested that there exists two models of theory organization; he stressed that, beyond thè "rational" model, an "empirical" model exists. Subsequently, for illustrating these two different models Lazare Carnot devoted two pages of each of his two books on mechanics; in particular, he claimed to develop his theory in an "empirical" way. In past century both Henri Poincaré and independently Albert Einstein again suggested two kinds of theory organization which were similar to the previous ones. For the bibliography, see (Drago 2012)

All in the above illustrates a basic dichotomy in the foundation of science; it concerns the philosophical notion of the organization of a theory and at the same time the corresponding formal theory of logic: either a problem-based organization (PO) with IL, or an axiomatic organization (AO) with CL.

There exists one more dichotomy which was long time debated by Intuitionists and Costructivists in opposition to Formalists about the two kinds of infinity, either the potential infinity (PI) or the actual infinity (AI). After some decades of this debate, this dichotomy was translated in formal terms of mathematics, obtaining a constructive mathematics (Bishop 1967) in opposition to the classical one.

In total, we have two dichotomies as the foundations of a scientific theory. Being AI and PI mutually incompatible, as well as AO and PO are (for philosophical and formal terms), two theories which are based on different choices are called incommensurable.

    **4. Kolmogorov's 1932 paper and a structural foundation of IL**

In the light of the two dichotomies Hilbert's program is easily recognized as founded on the pair of choices on these dichotomies AI&AO. By having well-defined this pair of choices Hilbert started the work of suggesting a structural foundation of the entire logic. Yet, he believed that this pair represented the unique way of founding logic (and science in general).

Instead the intuitionist program was founded on subjective notions, which however implicitly represented the alternative pair, PI&PO, i.e. the opposite pair of choices to Hilbert's. By going beyond Luitzen Brouwer's personal, subjective thinking mathematics and logic, in 1932 Andrei Kolmogorov suggested a foundation of IL as based on an objective pragmatic problem, the resolutions of mathematical problems through a calculus of them. The author illustrated this theory in a way which closely corresponds to the pair of choices of the intuitionist program. By a rational reconstruction of Kolmogorov's system according to both the choice for the PO model and the choice for the PI, I obtained one more instance of a rigorously intuitionist reasoning within a scientific theory (the theory of IL); and also I obtained a structural foundation of IL, to be put on par of an AO theory of it. In such a way logic is founded in a structural way, going beyond the reference to the objective meanings of some crucial notions (eg inference, negation, etc.), problems. (Drago 2020).

Hence, there exist two separate formal formulations of 'general logic', one is CL and the other is IL. In other words, the work of founding a logic is essentially of a pluralist nature.

**5. Founding logic by means of a structural problem: What is Negation?**

However, Kolmogorov's foundation of IL is based on a problem of an objective nature: how define a calculus of the resolutions of mathematical problems. Decades later, Michael Dummett (1977) has suggested a new foundation: it is based on the problem of how verify intuitionist propositions through proofs:

> From an intuitionist standpoint…. an understanding of a mathematical statement consists in the capacity to recognize a proof when presented with one; and a proof of such a statement can consist only in the existence of such a proof. (Dummett 1977, p. 6)

Again this problem is of an objective nature. A foundation which is based on a problem concerning a basic feature or element of IL is not still known.[6] In the following such a kind of foundation will be suggested for what Kolmogorov (1924/25, p. 418) called "general logic', i.e. informal logic including both CL and IL. By clarifying what is a negation within general logic we will obtain the foundations of both formal kinds of logic, CL and IL.

Being impossible an axiomatic (AO) within the informal context of general logic only a PO of this theory is possible.

---

[6] The above foundation of IL does not follow Brouwer's views; he considered logic as the regularities of mathematical arguing; in particular, he represented the absurd through a mathematical equality: $0=1$. Here, IL will be rather founded as a self-reliant theory.

The first step of the theoretical development of a PO theory consists in stating the universal problem upon which organize the theory of general logic. Let us take the problem "What is a negation?", which refers to a basic logical notion. Surely, the negation is a unary operation on either a term or a proposition; its result is a corresponding (negative) term or proposition which are different from the previous ones.

The second step consists in recognizing as its background the common logical knowledge; exactly what the PO theory of Kolmogorov's 1932 paper did through the two lists A and B concerning the calculus of problems (pp. 330-331) [7]. We can reiterate Kolmogorov's lists, apart 4.1 and 4.1.1 of A, both concerning negation.

Let us remark that one can reiterate this operation of negation an unlimited number of times; we obtain an infinite number of propositions. That gives the problem of how attribute meanings to all these propositions. A simple way to approach such a uncontrollable situation is to focus the attention on the doubly negated proposition and ask the simple question: Is this doubly negated proposition the same (i.e. not different) from the affirmative proposition, or not?

The answer 'same' gives the DNL and moreover leads to found CL through it. The other answer, 'different', leads to found NCL. So, the first decision of general logic is to establish not the truth or the falsity of all well-formed propositions (as CL does) but the sameness or not of a doubly negated proposition with its corresponding affirmative proposition.[8] In this way the infinite number of $n$-times negated propositions collapses into only three propositions: affirmative, negative and doubly negated. All that proves that the general notion of double negation precedes a specific kind of logic; it leads to a *choice* on the kind of logic. Conversely, its meaning is determined by the choice performed: either an affirmation, or a surmise[9].

As a consequence, also a negation takes two specific and mutually incompatible meanings. Within CL, since two negations are equivalent to the corresponding affirmation (and all other $n$-times negated propositions are equivalent to either the affirmative corresponding proposition or the negative proposition), negation plays the role of representing the entire logical world outside the corresponding affirmation; in other words, it is complementary in a clear-cut way to the affirmation. That implies that the opposition between the negation and the affirmation is of an antagonistic, exclusive nature. This fact leads to establish both the bivalence principle on all proposition and in general the non-contradiction principle: "Either A or not A".

Instead within IL, negation cannot oppose to the affirmation because it represents an open situation[10]; the proposition: "Acquitted for insufficiency of evidence of guilty" does not means a correct behavior of the defendant and a misconduct either.

In general one can attribute to the intuitionism the following principle:

Any [affirmative] proposition that is not without content should refer to one or more completely determine states of affairs accessible to our experience. (Kolmogorov 1932, p. 332),

Or, more shortly: "Any intuitionist affirmative proposition is supported by full evidence". That does not occurs in the case of a negation:

*For any universal proposition it is in general meaningless to consider its negation as a determinate proposition.* (Kolmogorov 1932, p. 333)

---

[7] It is also what Lobachevsky did by carelessly listing fifteen preliminary propositions (the propositions no.s 4, 7, 9 are deducible form the previous ones) to his most relevant illustration of non-Euclidean geometry. and

[8] However, also within IL a sameness is established, yet about the next number of negations of a proposition, three: ⌐⌐⌐ a = ⌐ a. Kolmogorov suggested to derive this equality from (a → b) → (⌐ b → ⌐ a), by replacing b with ⌐⌐ a (being of course a→ ⌐⌐ a).

[9] In complete isolation, Cusanus (1401-1464) anticipated some laws of IL (of terms) by reasoning about theological subjects (Drago 2017b). He stated that "… a surmise is a positive assertion that partakes—with a degree of otherness—of truth as it is." (Cusanus 1441-1442, I, 11, no. 57, p. 190) The word "partakes" summarizes a DNP because it is equivalent to: "It is not the case that it is not…".

[10] Notice that "Negation does not oppose to affirmation" was stated the first time by Cusanus. Remarkably, he also stated: "the negation does not exists before the not-other." (Cusanus 1462, chp. 4, p. 1113, no. 11); and not other is not the same (*ibidem*, chp. 5, p. 1117, no. 19; p. 1164, no. 123)

The same hold true for affirmation, because double negation holds true exactly when the evidence of the affirmation is insufficient; that means that also affirmation partakes absurdity (as it is apparent if one states: "Defendant's behavior was correct" in the case of insufficient evidence for his guilty).

Therefore, the definition of the logical operation of negation depends from the different kind of organization of the logical theory, either the AO governed by CL which leaves no room to double negation, or the PO governed by NCL, where double negation is more relevant than affirmation. Owing to this radical variation in meaning of the shared basic notion of negation, CL has to be considered as incommensurable with IL and other kinds of logic[11]. All that gives reason why scholars wanting to define negation by considering no more than internal operations of a the specific kind of logic met with insurmountable difficulties.

The next step of the specific development of a PO theory is an AAA. Here, it is motivated by the naturally next question, i.e. whether negation is connected to arguments resulting into the absurd; being impossible a direct proof solving this question, one has to apply an AAA. It starts by negating the thesis; i.e. by putting as hypothesis of a subsequent argumentation the following proposition "Negation <u>does not</u> partake absurdity". Yet, let us remark that affirmation cannot of course play the role of separation between the logical discourse and absurdity; double negation either, because within the specific logic of general logic, CL, the double negation version equates the corresponding affirmation which is separated from absurdity. Among the three versions of a proposition, only the negative one remains. If also negation could not imply the absurdity, then the separation of a logical discourse from absurdity is lost and we cannot avoid that this discourse includes absurd propositions; that means that it is essentially confuse and hence pointless. Hence, we have to deny the denied thesis of previous AAA; we obtain the DNP: "It is <u>not</u> the case that a negation <u>does not</u> partake of absurdity".

Within CL the application of DNL to this result of the theoretical development gives the proposition: "Negation partakes absurdity" as a valid proposition in general. Since by definition affirmation (and hence the double negation too) is absurdity-free, negation includes the entire absurdity. In other terms, given the clear-cut separation between affirmation and negation, the partaking of absurdity by negation implies its full relation with absurdity. Hence, the above conclusion is changed into the following proposition: "Negation is absurdity". Within the AO theory of logic, where all is mutually derived, it is more appropriate to represent this relation (a copula) with absurdity through an implication – i.e. one of the two mutual implications (direct and inverse) summarized by a copula: "Negation implies absurdity". In this case the implication proof exists, also if it is intended in an idealistic way, because in CL idealistic existence does not lead to contradiction.

In sum, CL reduces the situation of the three versions of a proposition to a sharply defined and very simple situation; the doubly negated proposition is actually cancelled because it is the same affirmative proposition; whereas the negated proposition is defined as the mirror black image of the affirmative proposition. At last, the logical arguing only has to address attention to, and entirely operate on the affirmative proposition.

Moreover, the usual illustration of CL translates this reference of a negation to an unknown totality of absurdity to a local opposition; it calls through a specific name, "false", what, by definition of negation, is "not true". CL then recuperates the general law *ex falso [absurdum et igitur] quodlibet*; where the intermediate words about absurdity are deliberately cancelled. In such a way within CL the absurdity becomes a consequence of the definition of false; this move exorcizes absurdity as a possible occurrence along the entire development of the discourse; absurdity only occurs after having met with a false proposition as a local step and a second-step event. The separation from absurdity has been therefore established in premise to the entire arguing

---

[11] Actually, one can conclude the existence of this incommensurability phenomenon since the previous remarks on double negation; which is an affirmation within CL and a surmise within NCL.

Moreover, the arguing is governed by two principles, bivalence and non-contradiction, which can decide in a clear-cut way the truth of the results of all derivations. No more simple logical arguing is possible. This fact justifies the dominant role played by CL and the kind of organization it governs, AO, along two millennia and more.

Let us come back to previous AAA's conclusion for establishing a structural foundation of IL according to the choice PO: "It is <u>not</u> the case that a negation <u>does not</u> partake absurdity" actually it is two times a DNP, because "partakes" means "share in part only", i.e. "it is <u>not</u> the case that it is <u>not</u> …".[12] Let us recall that within IL the DNL fails; it cannot be applied to this conclusion. Rather, according to the model of a PO theory, we can apply PSR; which, by changing at once the double negation of the proposition and the double negation of "partakes", gives the same proposition as the above, valid one within CL: "Negation is absurdity". However its present context is not that of a certainty, but that of an open, in progress situation, because the application of PSR is not a logical inference bringing with it a compelling logical force; rather, it is a translation between two kinds of logic; it merely appeals to a subjective belief on the completeness of the human reason; as such, it works as a suggestion of hypotheses, not certainties.

That means that the resulting theory of logic (IL) represents not a pyramid of infinite in number, static truths, like those which are obtained by applying CL within an AO theory; but an investigative process aimed at testing its final results with reality.[13]

### 6. The *ad absurdum* proof implied by a negation: Possibilism or actualism?

This difference between CL and IL about a negation is a subtle one, it concerns the method – either DNL or PSR - for obtaining an apparently same conclusion. The first method considers a negation as implying absurdity, i.e. as suggesting the existence of a proof of its absurdity; whereas the latter method, since achieves the common affirmative conclusion by means of a change of the kind of logic supported by an appeal to the human rationality, suggests the hypothesis of the existence of a proof of absurdity; which is not always manifest. In the first case the proof is actual, in the second case the proof is no more than possible.

> Thus the negation of a is transformed into an *existential proposition:* There exists a chain of logical inferences that, under the assumption of the correctness of a, leads to a contradiction. (Kolmogorov 1932, p. 332)

This subtle question actually represents the main basic problem of present philosophy of IL. Brouwer suggested that a negative proposition exactly implies, through a specific construction the absurdity.

> Brouwer [from 1923a-c]… expressed negation as reasoning that leads to an absurdity, or, briefly, as an absurdity.(Franchella 1994, p. 258).[14]

By remarking that this formalization implies to have in precedence obtained the notion of negation in the AAA or in the law of contradiction, Kolmogorov (1924/25, §§3-6, pp. 420-422) distrusted in it. He discussed how formalize a negation; yet his arguing is disputable.

This question remained unresolved also after a century long debate (Sundholm 1994; Raatikainen 2004; Raatikainen 2013). The second author concluded:

> We have examined the three basic choices there for the intuitionist theory of truth, the strict actualism, the liberalized actualism and possibilism, and found all them wanting. (Raatikainen 2004, p. 143)

I suggest that the question is unresolved because it was scrutinized with reference to connectives and proofs, but not the entire organization of a theory; it was scrutinized inside the AO, i.e. without taking in account the PO.

---

[12]  Being a cumbersome proposition, one may state it as the corresponding S4 modal proposition: "Negation may partake absurdity" ("Ex negatio quodlibet posse sequi").

[13]  Notice that this dichotomic view of the entire theories is not apparent within a single proposition or implication of a specific logic, apart from the doubly negated propositions. Once these propositions are underestimated or exorcized Cl maintains his monopolistic domain on logic.

[14]  The book (de Stigt 1990, pp. 238-270) offers a more detailed analysis of the subject.

Let us consider a clear instance of a PO theory, Lobachevsky's theory of parallel lines (1950). After fifteen preliminary propositions, the proposition no. 16 put the problem of how many are the parallels lines; "in the uncertainty" he decides to explore the case of two parallels lines. First, he suggests a new, more general definition of parallelism; then he states two theorems corroborating his definition and then he proves four theorems whose general conclusion is the following proposition:

> "The second assumption [of two parallel lines] can likewise be admitted <u>without</u> leading to any <u>contradiction</u> in the results..."[15]

(This proposition is a DNP since it is not equivalent to: "The second assumption is consistent"). This proposition represents a surmise on the possibility of correctly reasoning on the base of two parallel lines. In the subsequent propositions he actually argues on the base of the corresponding affirmative proposition; that means that he has applied the PSR to the previous DNP. Afterwards, from it he deductively derives fourteen theorems of his non-Euclidean geometry. (Bazhanov and Drago 2010)

Lobachevsky never considered as proven his new geometry. His search for an astronomical confirmation was unsuccessful. Moreover, he had no decisive mathematical proof of the truth of his geometry (only the translation of the Euclidean spherical trigonometry into the hyperbolic spherical trigonometry through the mere change of the angle α into the imaginary angle iα). For this reason he called his geometry "imaginary", a word to be meant in the same sense of the imaginary numbers; i.e. numbers which are not real, yet through them it is possible to obtain real results.

Actually, this is exactly the sense of a PO theory: the conclusion of the final AAA suggests a surmise; but a surmise, being a DNP, cannot be accurately tested with reality; therefore the application of PSR translates this conclusion into an affirmation and the logic into CL and the entire PO theory into AO, in order to test all them with reality. If the answer obtained is positive, the entire theory becomes effective, i.e. what was 'imaginary' is then real. The same for the proof implied by a negation; this proof is possible before the theory is changed by the PSR; afterwards, when the theory is governed by CL, it is an actual proof.

As an application, let us consider Natural Deduction. Let us recall that Natural deduction is not completely well formalized:

> There are no clear definitions of what is a proof theoretic formulation of a logic… and what is e.g. a Gentzen formulation. (Gabbay 2014, p. 44).
>
> The notion of a proof system is not well defined in the literature. There are some recognized methodologies such as 'Gentzen formulations', 'tableaux', 'Hilbert style axiomatic systems', but these are not sharply defined. (Gabbay 2014, p. 45).

One can see a reason of the insufficiency of Natural deduction in its ill-definition of a negation; it is defined as the classical one, i.e. a negation implies absurdity in an actual way. That is, by going against the assumptions of IL, ND assumes in all cases that negation implies the existence of a proof of absurdity, also when this proof is not known. Only this choice allows current Natural Deduction to attribute the same introduction and elimination rules to the negation in both kinds of logic.

In sum, present ND ignores that the difference between classical logic and intuitionist logic is not only a logical law, i.e. the DNL, but also the definition of negation, which cannot be the same for the two kinds of logic.

## 8. Conclusion

In a well-known paper concerning the difference between CL and other kinds of NCL, Prawitz recalled:

> A basic tenet of intuitionism is that classical logic contains some invalid forms of reasoning and consequently has to be rejected and, at least within mathematics, replaced by intuitionistic logic. In discussions of intuitionistic logic the question of the validity of this claim is often evaded, and instead intuitionistic logic is justified as being of interest from same special point of view which does

---

[15] By passing, we remark that in all his writings about non-Euclidean geometries Lobachevsky never stated "There exist two parallel lines"; instead he always wrote a DNS, like the previous one.

not necessarily repudiate the canons of classical logic but allows the peaceful coexistence of the two systems. Howcver, for anybody concerned with logic, the question whether the intuitionistic attack on classical logic is justified must be a vital issue. (Prawitz 1977, p. 2)

The present paper suggests that the attack of intuitionism was inappropriate; exclusivists are the choices, not the theories. A co-existence is therefore the solution, but not in technical terms of indifferently adopting all technical tools, but distinguishing the theories according to the basic choices on the two dichotomies; in particular by considering the radical variations of meaning given by the two basic choices, AO and PO, on the kind of organization of a theory.

**Bibiliography**